\input amstex
\documentstyle{amsppt}
\document

\def \NN{\Bbb N}
\def \ZZ{\Bbb Z}
\def \RR{\Bbb R}
\def \CC{\Bbb C}

\hcorrection{2cm}
\vcorrection{2cm}

\centerline{\bf EXTENSION AUX CYCLES SINGULIERS }\par
\vskip 3mm
\centerline{\bf  DU THEOREME DE KHOVANSKI-VARCHENKO.}\par

\vskip 7mm

\centerline{ par Abderaouf Mourtada}

\vskip 2cm

\heading Université de Bourgogne, I.M.B.\\
        U.M.R. 5584 du C.N.R.S., U.F.R. des Sciences et Techniques\\
         9, avenue Alain Savary, B.P. 47 870, 21078 Dijon Cedex.\\
      \rm E-mail: mourtada\@u-bourgogne.fr  \endheading

\vskip 3cm

{\bf Abstract.}{\sl  Let $\omega=dH$ be a hamiltonian 1-form in the real plane, of degre $d$. In [K][V], Khovanski and Varchenko proved that for any algebraic unfolding $\omega_\nu$ of $\omega$, of degre $d'$, with non-vanishing Abelian integrals along real cycles of $\omega$, the number of limit cycles of $\omega_\nu$, which born from regular real cycles of $\omega$,is bounded by some function of the degres $d$ and $d'$. In this paper, we extend this result to singular real cycles (polycycles), assuming that $H$ is a Morse function on $\CC^2$. We deduce in particular the following result: if $d=d'$ and $H$ generic at infinity, then the number of limit cycles of $\omega_\nu$ in the real plane, is bounded by a function of the degre $d$.}\par

\vskip 2cm

\noindent{\bf Introduction.}

\vskip 3mm

 Soit $H:\RR^2\to\RR$ un polyn\^ome de degr\'e $d+1$, dont les points critiques de son complexifié sont de Morse (les valeurs critiques correspondantes ne sont pas forc\'ement distinctes deux \`a deux). Soit $H_{d+1}$ le bloc homog\`ene de degr\'e $d+1$ dans $H$. On suppose que\par

\vskip 3mm

\noindent ($*$) $H_{d+1}$ est un produit de facteurs $\CC$-lin\'eaires, deux \`a deux distincts.\par

\vskip 3mm

 Soit $\nu=(\epsilon,v)\in\RR^{(d+1)(d+2)}$ et soit $\omega_\nu=dH+\epsilon\eta_v$ un d\'eploiement alg\'ebrique de degr\'e $d$. Soit $({\Cal C}_i)$ la famille des couronnes de cycles r\'eels de la fibration $H$ et soit $I_i(\eta_0)$ l'int\'egrale de la 1-forme $\eta_0$ le long des cycles de la couronne ${\Cal C}_i$. On suppose que\par

\vskip 3mm

\noindent ($**$) $I_i(\eta_0)\not\equiv 0$ pour tout $i$.\par

\vskip 3mm

\noindent Le but de ce travail est d'\'etablir le r\'esultat suivant

\vskip 3mm

\proclaim{Th\'eor\`eme} Il existe un entier $N(d)$ ne d\'ependant que du degr\'e $d$ tel que, pour $\nu$ suffisament petit, la 1-forme $\omega_\nu$ a au plus $N(d)$ cycles limites (compt\'es avec multiplicit\'e) dans le plan r\'eel.
\endproclaim

\vskip 3mm

 Le nombre de ces couronnes ${\Cal C}_i$ est major\'e par une fonction du degr\'e ($<4d^2$). Le th\'eor\`eme se d\'eduit donc de la

\vskip 3mm

\proclaim{Proposition} Il existe un entier $n(d)$ tel que pour toute couronne ${\Cal C}$ de cycles r\'eels de $H$, il existe un voisinage ${\Cal V}_{H,\eta_0}$ de $\overline{{\Cal C}}$ dans le disque de Poincar\'e tel que, pour $\nu$ suffisament petit, $\omega_\nu$ poss\`ede au plus $n(d)$ cycles limites (compt\'es avec multiplicit\'e) dans ${\Cal V}_{H,\eta_0}$.
\endproclaim

\vskip 3mm

 Soit $\{\delta(t);\ t\in ]t_1,t_2[\}$ la famille des cycles r\'eels constituant ${\Cal C}$, le cycle $\delta(t)$ est une composante connexe de la fibre r\'eelle $\{H=t\}$; et les r\'eels $t_j$ sont des valeurs critiques de $H$ (si $\ne\infty$). Soit $I(t)=I_{\delta(t)}(\eta_0)$ l'nt\'egrale ab\'elienne de $\eta_0$ sur le cycle $\delta(t)$. Dans [K][V], Khovanski et Varchenko ont montr\'e que le nombre et la multiplicit\'e de z\'eros de l'int\'egrale $I$ sur l'intervalle $]t_1,t_2[$, sont major\'e par une fonction du degr\'e $n_1(d)$. Donc, par le classique lemme de perturbation [AL], le nombre et la multiplicit\'e des cycles limites de $\omega_\nu$ qui naissent \`a partir des cycles $\delta$, sont major\'es par $n_1(d)$. Le bord de la couronne ${\Cal C}$ dans le disque de Poincar\'e est l'union de deux composantes connexes. Par l'hypothèse $(*)$, chacune de ces composantes est, ou bien un cycle singulier $\Gamma_k$, \`a $k$ singularit\'es (auquel cas $\Gamma_k\subset\{H=t_j\}$), ou bien un cycle r\'egulier $\Gamma_0$ sur l'\'equateur (auquel cas $t_j=\infty$). Si $t_j\ne\infty$, il est connu (cf. [AV] par exemple), que l'int\'egrale Abélienne $I$ admet un d\'eveloppement asymptotique logarithmique au voisinage de $t_j$: $I(t)=\sum_{n,m\leq n} a_{n,m}(t-t_j)^n(\log (t-t_j))^m$. On appelle {\bf multiplicit\'e alg\'ebrique} de $I$ en $t_j$, et on la note $ma(I,t_j)$ le plus petit entier $n$ tel que $a_{n,m}\ne 0$ pour un certain $m$. Elle coincide avec la multiplicit\'e classique si $I$ est analytique au voisinage de $t_j$. Si $t_j=\infty$, on montre dans le §2, que la fonction $J(\tau)=\tau^{d+2}I(\tau^{-(d+1)})$ est analytique au voisinage de 0. Dans ce cas, on note $ma(I,\infty)=ma(J,0)$. La proposition est alors une cons\'equence des deux lemmes suivants

\vskip 3mm

\proclaim{Lemme 1} Il existe un entier $M(d)$ tel que pour tout $t\in[t_1,t_2]$, $ma(I,t)\leq M(d)$.
\endproclaim

\vskip 3mm

\proclaim{Lemme 2} Il existe un entier $n_2(d,ma(I,t_j))$ et un voisinage $V_{H,\eta_0}$ de $\Gamma_k$ dans le disque de Poincar\'e tels que, pour $\nu$ suffisament petit, $\omega_\nu$ a au plus $n_2$ cycles limites (compt\'es avec multiplicit\'e) dans $V_{H,\eta_0}$.
\endproclaim

\vskip 3mm

  La référence la plus générale sur les résultats concernant les intégrales Abéliennes est [R]. Dans [G3], Gavrilov montre un résultat globale (du type du théorème ci-dessus), dans le cas $d=3$, et sans l'hypothèse $(**)$. Dans ce cas, la monodromie de la fibration $H$ est {\bf transitive}, et les seuls cycles singuliers rencontrés sont du type $\Gamma_1$. La preuve de [G3] s'appuie sur un théorème de Roussarie [Ro], qui est un cas particulier du théorème principal 2.1 ci-dessous. Ce théorème est basé sur le théorème IVB1 de [Mo]. Une généralisation de la situation dans [G3], est la suivante: on suppose que la fibration de Morse $H$ satisfait à l'hypothèse $(*)$, et qu'elle est de monodromie transitive (le groupe fondamental de $\CC\setminus\{\text{les valeurs critiques de}\ H\}$ agit transitivement sur le groupe d'homologie de la fibre générique de $H$). Dans ce cas, on peut montrer un théorème global (où le voisinage ne dépend que de $H$), en utilisant les résultats de [B], et le théorème 2.1. Dans le cas général (monodromie non transitive), il est raisonnable d'étudier d'abord les déploiements de couronnes de cycles réguliers, dans l'esprit des derniers travaux de Gavrilov ([G4]...). La jonction vers les cycles singuliers se fera via le théorème IVC1 de [Mo].\par
  
\vskip 3mm

\noindent{\bf 1. D\'emonstration du lemme 1.}\par

\vskip 3mm

 Soit $H:\CC^2\to\CC$ une fonction polynomiale de degr\'e $d+1$ dont les points critiques $(m_{i,j})_{j=1,\ldots,p\ i=1,\ldots,\ell_j}$ sur $\CC^2$ sont de Morse et qui satisfait \`a l'hypoth\`ese $(*)$ (donc par le th\'eor\`eme de Bezout $\sum_{j=1,\ldots,p} \ell_j=d^2$). Notons $T=\{t_1,\ldots,t_p\}$ l'ensemble des valeurs critiques. Soit $V_t=H^{-1}(t)$ et $\overline{V}_t\subset \CC P^2$ sa clôture projective. Elle coupe transversalement la ligne \`a l'infini $\CC P^1$ en $d+1$ points ind\'ependants de $t$. Pour $t\in \CC\setminus T$, la fibre r\'eguli\`ere $V_t$ est {\bf d'homologie \'evanescente et born\'ee} de dimension $d^2$ sur $\ZZ$ ([AV], [I1], [I2],). Pour $j\in\{1,\ldots,p\}$, soit $\gamma_j\subset\CC$ un petit lacet autour de la valeur critique $t_j$ et $h_j:H_1(V_t,\ZZ)\to H_1(V_t,\ZZ)$ l'op\'erateur de monodromie "classique" correspondant ($t\in\gamma_j$). C'est un produit d'op\'erateurs de monodromie de Picard-Lefschetz ([AV], [G1], [G2]): en effet, soit $H_\lambda:\CC^2\to\CC$ une famille \`a un param\`etre $\lambda\in\CC$ de fonctions polynomiales de degr\'e $d+1$ telle que $H_0=H$ et pour $\lambda\ne 0$ suffisament petit, les valeurs critiques $(t_{i,j}(\lambda))$ de $H_\lambda$ soient distinctes deux \`a deux. Or, pour $t\in\gamma_j$, les fibres $V_t$ et $V_{\lambda,t}=H_\lambda^{-1}(t)$ sont isomorphes, et l'op\'erateur de monodromie de $H_\lambda$ correspondant au lacet $\gamma_j$ est le produit des op\'erateurs de monodromie de Picard-Lefschetz $h_{i,j}$ associ\'es aux valeurs critiques $t_{i,j}$. Soit $\Delta_{i,j}\in H_1(V_t,\ZZ)$ le cycle \'evanescent au point critique $m_{i,j}$ et $c\in H_1(V_t,\ZZ)$. Un calcul direct donne

$$
 h_j(c)=c+\sum_{i=1}^{\ell_j} a_i(c)\Delta_{i,j}
\tag 1
$$

\noindent les entiers $a_i(c)\in\ZZ$ d\'ependent des indices d'intersection des cycles $c$ et $\Delta_{i,j}$.\par

\vskip 3mm

 A la fibration globale de Milnor $H$ au dessus de $\CC\setminus T$, on associe la fibration homologique globale de Milnor $E$ de m\^eme base et dont les fibres $E_t$ sont les espaces homologiques $H_1(V_t,\CC)$. Soit $\eta$ une 1-forme alg\'ebrique complexe sur $\CC^2$ de degr\'e$\leq d$. Si $\delta$ est une section locale constante de $E$, on note $I_\delta(t)$ l'int\'egrale Ab\'elienne de $\eta$ sur le cycle $\delta(t)\subset V_t$. Il est connu ([AV]) que chaque branche de $I_\delta$ est analytique au voisinage de chaque point de $\CC\setminus T$. L'homologie de $H$ \'etant born\'ee, on a pour tout $j\in\{1,\ldots,p\}$

$$
 I_\delta(t)=O(1)
\tag 2
$$

\noindent sur le germe de tout secteur $S_j$ bas\'e au point $t_j$. Et, d'apr\`es [Ma], [Y], on a

$$
 I_\delta(t)=O(t)
\tag 3
$$

\noindent sur le germe de tout secteur $S_\infty$ bas\'e au point $\infty\in\CC P^1$. Dans [Y], on trouve aussi une estimation moins fine, mais établie dans le cas g\'en\'eral o\`u $H$ ne satisfait pas \`a l'hypoth\`ese $(*)$. Sa preuve est bas\'ee sur des notions \'el\'ementaires de th\'eorie de l'\'elimination pour la localisation des points de ramification associ\'es aux courbes alg\'ebriques affines $V_t$.\par

\vskip 3mm

 La ramification de $I_\delta$ au point $t_j$ s'obtient gr\^ace \`a la monodromie $h_j$ donn\'ee par (1): par la d\'efinition des cycles \'evanescents $\Delta_{i,j}$ aux points critiques de Morse $m_{i,j}$, les int\'egrales $I_{\Delta_{i,j}}$ sont analytiques aux voisinages de $t_j$ et tendent vers 0 quand $t$ tend vers $t_j$ ([AV]). Soit

$$
 J_{\delta}(t)=\frac{1}{2i\pi}(\sum_{i=1,\ldots,\ell_j} a_i(\delta)I_{\Delta_{i,j}}(t))\log (t-t_j)
$$

\noindent en appliquant la monodromie $h_j$ au cycle $\delta(t)$, on obtient que la fonction

$$
 f_\delta(t)=I_\delta(t)-J_\delta(t)
\tag 4
$$

\noindent est uniforme sur un voisinage point\'e de $t_j$, et par (2) elle est analytique au voisinage de $t_j$.\par

\vskip 3mm

 Soit $\{\delta_1,\ldots,\delta_{d^2}\}$ une base de sections locales constantes de $E$ et soit $W(t)$ la matrice wronskienne des int\'egrales $I_{\delta_j}(t)$ de rang $\ell$. Soit $w(t)$ un wronskien d'ordre $\ell$ non identiquement nul. Par une d\'emarche classique utilisant (2), (3) et (4) (cf. [I1], [I2], [I3], [Ma], [Y]), $w$ est une fonction rationnelle dont les p\^oles appartiennent \`a $T$ et dont les degr\'es du num\'erateur et du d\'enominateur sont major\'es par une fonction du degr\'e $d$. Ainsi, si $\delta$ est une section locale constante de $E$ telle que $I_\delta\not\equiv 0$, la multiplicit\'e alg\'ebrique $ma(I_\delta,t)$ en tout point $t\in\CC$, est major\'ee par une fonction du degr\'e $d$. Soit $M$ le plus grand des entiers $m$ tels que la fonction $t^{m-1}I_\delta(t)$ soit born\'ee sur un certain secteur $S_\infty$. Le m\^eme type de raisonnement appliqu\'e \`a $w(t)$ montre que $M$ est major\'ee par une fonction du degr\'e $d$, et ceci finit la preuve du lemme.\qed

\vskip 3mm

\noindent{\bf 2. D\'emonstration du lemme 2.}

\vskip 3mm

 Soient $X_0$ et $X_\nu$ les champs de vecteurs de $\RR^2$ associ\'es aux 1-formes $dH$ et $\omega_\nu$. Commen\c cons la preuve dans le cas d'un cycle r\'egulier \`a l'infini. Soit ${\Cal C}=\{\delta(t);\ t\in]t_1,+\infty[\}$ une couronne de cycles r\'eels $\delta(t)$ de $H$. Soit $\overline{X}_0$ et $\overline{X}_\nu$ les prolongements des champs $X_0$ et $X_\nu$ sur $\CC P^2$. Le champ $\overline{X}_0$ a un cycle r\'egulier $\Gamma_0=\RR P^1\subset \CC P^1$ qui ne porte pas de l'holonomie. Soit $(u=1/x,v=y/x)$ une carte sur $\CC P^1\subset\CC P^2$ et $\sigma$ le germe en 0 de la transversale $\{v=0\}$. Le champ $X_\nu$ \'etant de degr\'e $d$, son prolongement $\overline{X}_\nu$ admet une application de retour $p_{1,\nu}$ sur $\sigma$ qui est analytique dans les coordonn\'ees $(u,\nu)$

$$
 p_{1,\nu}(u)=u+\epsilon(K(u)+O(\nu))
\tag 5
$$

\noindent D'un autre c\^ot\'e, consid\'erons la coordonn\'ee r\'eelle $\tau=t^{-1/(d+1)}$ sur la semi-transver\-sale $\sigma^+=\RR^{+*}\cap\sigma$. L'application de retour $p_{2,\nu}$ sur $\sigma^+$ s'obtient en int\'egrant la 1-forme $d(H^{-1/(d+1)})$ le long des orbites du champ $X_\nu$. Un calcul direct donne

$$
 p_{2,\nu}(\tau)=\tau+\epsilon(\tau^{d+2}I(\tau^{-(d+1)})+O(\nu))
\tag 6
$$

\noindent o\`u $I(t)=I_{\delta(t)}(\eta_0)$ est l'int\'egrale de la 1-forme $\eta_0$ le long du cycle $\delta(t)$. Maintenant, la relation $\tau^{-(d+1)}=H(1/u,0)$ et l'hypoth\`ese $(*)$ montrent qu'il existe un diff\'eomorphisme $g$ analytique en 0 tel que $\tau=g(u)$. Par cons\'equent, en conjuguant (5) \`a (6) par $g$, on obtient que la fonction $J(\tau)=\tau^{d+2}I(\tau^{-(d+1)})$ est analytique en 0, et ceci conclut la preuve dans ce cas.\par

\vskip 3mm

 Soit maintenant ${\Cal C}=\{\delta(t);\ t\in]0,t_2[\}$ une couronne de cycles r\'eels de $H$ et $\Gamma_k=\partial{\Cal C}\cap\{H=0\}$ un cycle singulier \`a $k$ singularit\'e; il est compact d'après l'hypothèse $(*)$. S'il est r\'eduit \`a un point (un centre), l'int\'egrale Ab\'elienne et l'application de retour correspondantes sont analytiques en 0. La preuve est alors similaire \`a celle donn\'ee ci-dessus. Sinon, ses singularit\'es sont des points de selle. Dans ce cas, la preuve est bas\'ee sur les idées g\'en\'erales d\'evelopp\'ees dans le travail [Mo].\par

\vskip 3mm

 Dans [Mo], il est \'etabli que le nombre et la multiplicit\'e des cycles limites de $\omega_\nu$ proches de $\Gamma_k$, sont uniform\'ement major\'es (cf. théorème 0, [Mo]). Il s'agit ici de montrer que ces majorants ne d\'ependent que du degr\'e $d$. On reprend les notations de [Mo]. Soient $(x_1,\ldots,x_k)$ des coordonn\'ees analytiques ad\'equates sur des transversales $\sigma_j$ \`a $\Gamma_k$. Soit $\lambda_j(\nu)$ le germe qui déploie la j-ièmme connexion de $\Gamma_k$. L'application de transition (ou application de Dulac) du $j$-i\`emme coin de $\Gamma_k$ pour la 1-forme $\omega_\nu$, s'\'ecrit 

$$
x_{j+1}=d_j(x_j,\nu)-\lambda_j(\nu)\quad\text{avec}\quad d_j(x_j,\nu)=x_j^{r_j}(1+D_j(x_j,\nu))
$$ 

\noindent o\`u le germe $D_j$ est induit par un \'el\'ement de l'alg\`ebre $QR{\Cal H}^{1,(1,d'=(d+1)(d+2))}$ (cf. VA[Mo]), et o\`u $r_j(\nu)=1+\mu_j(\nu)$ est le nombre caract\'eristique de la j-i\`emme singularit\'e de $\Gamma_k$. Si on note $g_j=d_j-x_{j+1}$ pour $j=1,\ldots,k-1$, et $f=d_k-x_1-\lambda_k$, alors l'application de retour de $\omega_\nu$ sur la transversale $\sigma_1$ (par exemple!) s'\'ecrit

$$
p_{1,\nu}(x_1)=x_1+f_{|\{g_1=\lambda_1,\ldots,g_{k-1}=\lambda_{k-1}\}}
\tag 7
$$

\noindent (cf. paragraphe IVC4, [Mo]). Dans la coordonn\'ee $t$ sur $\sigma_1$, cette application de retour s'obtient en int\'egrant la 1-forme $dH$ le long des orbites de $\omega_\nu$

$$
p_{2,\nu}(t)=t+\epsilon(I(t)+O(\nu))
\tag 8
$$

\noindent o\`u $I(t)$ est l'int\'egrale Ab\'elienne de $\eta_0$ sur les cycles $\delta(t)$. D'après [AV], c'est un \'el\'ement de l'alg\`ebre convergente $QR{\Cal H}_{cvg}^{1,0}$ (cf. IA[Mo]) . Cependant, le terme $O(\nu)$ n'est pas uniforme dans la variable $t$ au voisinage de 0. L'\'ecriture (7) de l'application de retour est donc la plus adapt\'ee pour \'etudier le probl\`eme, et elle s'interpr\`ete g\'eom\'etriquement de la fa\c con suivante: pla\c cons nous dans un voisinage $U$ de 0 dans $(\RR^{+*})^k\times\RR^{d'}$ de coordonn\'ees $(x=(x_j),\nu)$, et dans lequel les germes $g=(g_j)_{j=1,\ldots,k-1}$ et $\lambda=(\lambda_j)_{j=1,\ldots,k-1}$ sont réalisés. La forme

$$
d\nu_1\wedge\cdots\wedge d\nu_{d'}\wedge dg_1\wedge\cdots\wedge dg_{k-1}
$$

\noindent ne s'annule pas sur $U$ (en réduisant $U$ si nécessaire). Soit $\chi$ la d\'erivation sur $(U,0)$ d'int\'egrales premi\`eres $\nu_1,\ldots,\nu_{d'},g_1,\ldots,g_{k-1}$, telle que $\chi x_1=\prod x_j$. Soit $U_0=\{(x,\nu)\in U;\ g(x,\nu)=\lambda(\nu)\}$), c'est une sous-variété analytique de $U$, invariante par le flot de $\chi$ dans $U$. Le nombre de points fixes de l'application de retour $p_{1,\nu}$ est \'egale au nombre de composantes connexes de l'intersection de la fibre $f_{|U_0}=0$, et de l'orbite de $\chi$  $\gamma_\nu$ correspondante aux valeurs $(\nu,\lambda(\nu))$ des intégrales premières de $\chi$. La multiplicit\'e de ces points fixes est major\'ee par l'indice de noeth\'erianit\'e de l'id\'eal diff\'erentiel restriction $I_{\chi,f|U_0}$.

\vskip 3mm

  Le cadre général de cette situation est le suivant (cf IVB, [Mo]): soient $k,q_1,q_2\in\NN$ et soient $(x,\alpha=(\mu,\nu))$ des coordonnées analytiques locales sur $((\RR^{+*})^k\times\RR^{q_1}\times\RR^{q_2},0)$. Soit un ouvert $U\in ((\RR^{+*})^k\times\RR^{q_1}\times\RR^{q_2},0)$ et soit $\chi$ une d\'erivation d'Hilbert réalisée sur $U$ (cf. IVA, [Mo]), de dimension de non trivialit\'e $k-1$, et d'int\'egrales premi\`eres non triviales $g=(g_1,\ldots,g_{k-1})$.\par
  
\vskip 3mm

  Soient $\lambda=(\lambda_1,\ldots,\lambda_{k-1})$ les valeurs (coordonn\'ees) des int\'egrales premi\`eres $g$ et soit $\gamma$ l'orbite de $\chi$ dans $U$, qui adhère à 0, et le long de laquelle toutes les int\'egrales premi\`eres de $\chi$ sont nulles. C'est une orbite principale de $\chi$ dans $U$ (cf. IB, [Mo]). On peut supposer, sans perte de généralité, que $U$ est le saturé par le flot de $\chi$ (dans $U$!) d'une transversale analytique à $\gamma$, qu'on note $\sigma$ et qu'on munit des coordonnées analytiques $(\alpha,\lambda)$. Soit $\pi_\chi: U\to \sigma$ la projection intégrale (on note de la même façon son germe aux points 0 et $\sigma\cap\gamma$).Soit $W_0\subset\sigma$ un semi-analytique de $\RR\{\alpha,\lambda\}$, qui adhère à 0, et soit $U_0=\pi_\chi^{-1}(W_0)$.\par
 
\vskip 3mm

  Soit $f\in QR{\Cal H}^{k,q}$ réalisé sur $U$, et soit $J_{\chi,f,\gamma}\subset\RR\{\alpha,\lambda\}$ son id\'eal $\chi$-transverse  le long de $\gamma$, c'est la restriction de l'id\'eal diff\'erentiel $I_{\chi,f}$ \`a une transversale analytique \`a $\gamma$ (cf. IB, [Mo]). On suppose que $f$ satisfait \`a l'hypoth\`ese $(H\lambda)$: il existe un entier $N$ tel que $J_{\chi,f,\gamma}\supset {\Cal M}_\lambda^{N}$ par restriction \`a $W_0$. Soit $N_0$ le plus petit de ces entiers $N$. D'apr\`es le th\'eor\`eme IVB1 de [Mo], le degr\'e de la projection $\pi_\chi$ restreinte \`a la fibre nulle de $f_{|U_0}$, est fini. De plus, l'id\'eal diff\'erentiel $I_{\chi,f}$ est localement noeth\'erien sur $U_0$. Notons $ind(I_{\chi,f|U_0})$ la borne sup\'erieure de son indice de noeth\'erianit\'e sur un voisinage choisi de 0 dans $U_0$. Si $ma=ma(\chi,f)_0$ est la multiplicit\'e alg\'ebrique du couple $(\chi,f)$ en 0 (voir ci-dessous), alors ce th\'eor\`eme se pr\'ecise de la fa\c con suivante

\vskip 3mm

\proclaim{Théorème principal 2.1} Il existe des fonctions universelles $N(ma,N_0,k,q_1)$ et \par $L(ma,N_0,k,q_1)$ telles que

$$
 d°\pi_{\chi|Z(f)\cap (U_0,0)}\leq N\quad\text{et}\quad ind(I_{\chi,f|U_0})\leq L
$$

\endproclaim

\vskip 3mm

 Indiquons bri\`evement comment on d\'eduit le lemme 2 du théorème 2.1 (cf. §2.2): soit $(\chi,f)$ le couple d'Hilbert (voir ci-dessous) associ\'e au d\'eploiement $\omega_\nu$ au voisinage du cycle singulier  $\Gamma_k$. En conjuguant (7) \`a (8), on montre que $N_0=1$ et $ma(\chi,f)_0 =ma(I(t),0)$. Par le th\'eor\`eme de Bezout,  les entiers $k$ et $q_1$ sont major\'es par une fonction du degr\'e $d$.\par

\vskip 3mm

\noindent{\bf 2.1. D\'emonstration du théorème 2.1.}\par

\vskip 3mm

\noindent{\bf (a). Multiplicit\'e alg\'ebrique g\'en\'eralis\'ee et son uniformit\'e.}\par

\vskip 3mm

 Soit $f\in QR{\Cal H}_{cvg}^{1,q=(q_1,q_2+\ell)}(x,\alpha=(\mu,\nu),u)$ et soit 

$$
\chi=x\frac{\partial}{\partial x}-\sum_{j=1}^\ell s_j(\mu) u_j\frac{\partial}{\partial u_j}
$$ 

\noindent avec $s_j(0)=1$. L'orbite $\gamma_0=\{ (\alpha,u)=0 \}$ est principale dans un ouvert $U\in (\RR^{+*}\times\RR^{|q|},0)$. On suppose que $f$ est réalisée sur $U$, continue sur $\overline{U}$. Il est montr\'e dans III[Mo] que $f$ admet une multiplicit\'e alg\'ebrique $ma(\chi,f)_0$ relativement \`a $\chi$ le long de $\gamma_0$: soient $\lambda=(\lambda_1,\ldots,\lambda_\ell)$ les valeurs (coordonnées) des intégrales premières $L_j=x^{s_j}u_j$.  Si $f=\sum_{p\in\ZZ} g_p$ est la s\'erie de $f$ en blocs $\chi$-homog\`enes et si on pose $f_n=\sum_{p\leq n} g_p$, la suite des id\'eaux $\chi$-transverses  $(J_{\chi,f_n,\gamma_0})_{n\geq n_0}$ de $\RR\{\alpha,\lambda\}$ est croissante, et son indice de stationnarit\'e (ind\'ependant de $n_0\in\ZZ^-$ suffisament grand en module) est la multiplicit\'e alg\'ebrique $ma(\chi,f)_0$. Or chaque orbite $\gamma_{m}=\{ (\mu,u)=0,\ \nu=\nu_m \}\subset U$ est aussi principale dans un voisinage dans $U$, du point $m=(0,0,\nu_m,0)$ du bord $B_{1,0}=\{ (x,\mu,u)=0\}\cap\overline{U}$. Et, par le même raisonnement que ci-dessus, le germe $f_{m}$ admet aussi une multiplicit\'e alg\'ebrique $ma(\chi_{m},f_{m})_{m}$ relativement \`a $\chi_{m}$ le long de $\gamma_{m}$ (en considérant les idéaux $\chi$-transverses dans l'anneau $\RR\{\mu,\nu-\nu_m,\lambda\}$). On note simplement $ma(\chi,f)_{m}$ cette multiplicité. La multiplicit\'e alg\'ebrique positive associée est $ma^+(\chi,f)_{m}=\max\{ma(\chi,f)_{m},0\}$. Un r\'esultat basique pour la suite est le

\vskip 3mm

\proclaim{Lemme 2.1.1} L'application $m\in B_{1,0}\mapsto ma^+(\chi,f)_{m}$ est semi-continue sup\'e-\par\noindent rieurement.
\endproclaim

\vskip 3mm

\noindent{\bf Preuve.} En $m=0$ par exemple. Par la définition de l'idéal $\chi$-transverse le long de $\gamma_0$, la multiplicit\'e alg\'ebrique $ma(\chi,f)_0$ est major\'ee par l'indice de stationnarit\'e de la suite croissantes des id\'eaux $(I_{\chi,f_n}(0))_{n\geq n_0}$ dans l'anneau $SB^{1,|q|}$. Notons $ma^+_0=ma^+(\chi,f)_0$. Le germe $f_{ma^+_0}$ a le même idéal $\chi$-transverse que $f$ et sa multiplicité algébrique positive est $ma^+_0$. Donc, par le théorème principal IIIA1[Mo], on a l'inclusion

$$
I_{\chi,f_{ma^+_0}}(0)\supset (x^{ma^+_0+\epsilon})\pi_\chi^*(J_{\chi,f,\gamma_0})\quad\text{pour tout }\ \epsilon>0
$$

\noindent où $\pi_\chi$ est (le germe en 0 de) la projection intégrale. Soit $h=f-f_{ma^+_0}$. On a $J_{\chi,h,\gamma_0}\subset J_{\chi,f,\gamma_0}$. Comme $f,h\in QR{\Cal H}_{cvg}^{1,q}$, la remarque IIIA1[Mo] donne

$$
(x^{n(f)}I_{\chi,h}(0)\subset \pi_\chi^*(J_{\chi,f,\gamma_0})
$$

\noindent et donc

$$
 (x^{n(f)+ma^+_0+\epsilon})I_{\chi,h}(0)\subset I_{\chi,f_{ma^+_0}}(0)
$$

\noindent Par conséquent, il existe $L\in\NN$ et $H_j\in SB^{1,|q|}$ tels que

$$
x^{n(f)+ma^+_0+\epsilon} h=\sum_{j=1}^L H_j \chi^j f_{ma^+_0}
$$

\noindent En germifiant cette égalité en tout point $m$ du bord, voisin de 0, on obtient

$$
 (x^{n(f)+ma^+_0+\epsilon})I_{\chi,h}(m)\subset I_{\chi,f_{ma^+_0}}(m)
$$

\noindent et par la d\'efinition de l'id\'eal $\chi$-transverse le long de $\gamma_m$, on a $ma^+(\chi,f)_{m}\leq ma^+_0$.\qed

\vskip 3mm

   Soit $p\in\NN$ et soient $(x,\rho,\alpha=(\mu,\nu))$ des coordonnées analytiques locales sur 
$((\RR^{+*})^{k+p}\times\RR^{q_1}\times\RR^{q_2},0)$. Soit $\chi$ une d\'erivation d'Hilbert, de dimension de non trivialit\'e $k-1$, et dont les int\'egrales premi\`eres non triviales sont 

$$
g_j(x,\rho,\alpha)=x_j^{r_j}(1+D_j(x_j,\rho,\alpha))-x_{j+1}\quad j=1,\ldots,k-1
$$

\noindent avec $r_j=1+\mu_j$ et $D_j\in QR{\Cal H}^{1+p,q)}$, $D_j(0,\rho,\alpha)\equiv 0$. Si $f\in QR{\Cal H}^{k+p,q}$, on dit que le couple $(\chi,f)$ est un {\bf couple d'Hilbert} de l'algèbre $QR{\Cal H}^{k+p,q}$. Ces  {\bf composantes} sont les germes $D_1,\ldots,D_{k-1}$ et $f$. En cas d'ambiguité, les "puissances" $r_j$ et la dimension de non trivialité $k-1$ seront précisées. Soit $(\chi,f)$ un couple d'Hilbert {\bf convergent}: ie. ses composantes appartiennent \`a l'alg\`ebre convergente $QR{\Cal H}^{k+p,q}_{cvg}(x,\rho,\alpha)$. On suppose que la dimension de non-trivialité du couple $(\chi,f)$ est $k-1$, et que les coordonnées $\rho$ sont des intégrales premières de $\chi$. Dans ce cas, on peut prendre $p=0$ pour simplifier, quitte à remplacer l'algèbre $QR{\Cal H}^{k+p,q}_{cvg}(x,\rho,\alpha)$ par une autre algèbre convergente $QR{\Cal H}_{cvg}^{k,q'}(x,\alpha')$ (en utilisant les fonctions élémentaires dans les coordonnées $\rho$ (cf. IA[Mo]). On suppose que ce couple est réalisé sur un ouvert $U\in ((\RR^{+*})^{k+p}\times\RR^{q_1}\times\RR^{q_2},0)$, continu sur $\overline{U}$.\par

\vskip 3mm

  Soit $(\pi,{\Cal N})$ le premier \'eclatement de $\chi$ de diviseur exceptionnel $\overline{{\Cal D}}=\overline{{\Cal D}(0)}$ (cf. IVA[Mo]), et soit $(\widetilde{\chi},\widetilde{f})$ le relev\'e de $(\chi,f)$. Dans les coordonn\'ees $(\rho,\alpha,u)$ au dessus de ${\Cal D}$, le champ relev\'e s'\'ecrit

$$
\widetilde{\chi}=\rho\frac{\partial}{\partial \rho}-\sum_{j=1}^{k-1} s_j(\mu) u_j\frac{\partial}{\partial u_j}
$$ 

\noindent avec $s_j(0)=1$ et $\widetilde{f}_a\in QR{\Cal H}_{cvg}^{1,.}(\rho,.)$ pour tout $a\in{\Cal D}$. On d\'eduit facilement du lemme 2.1.1

\vskip 3mm

\proclaim{D\'efinition et Lemme 2.1.2} La multiplicit\'e alg\'ebrique du couple d'Hilbert $(\chi,f)$ en chaque point $m=(0,0,\nu_m)\in B'_{k,0}=\{(x,\mu)=0\}\cap\overline{U}$, est la multiplicit\'e alg\'ebrique positive $ma^+(\widetilde{\chi}_{a_0},\widetilde{f}_{a_0})_{a_m}$ au point $a_m=(0,0,\nu_m,0)\in{\Cal D}$. On la note $ma(\chi,f)_{m}$. L'application $m\in B'_{k,0}\mapsto ma(\chi,f)_{m}$ est semi-continue sup\'erieurement.
\endproclaim

\vskip 3mm

 Soient $t$ un r\'eel$>0$ et $\psi_\mu$ une famille analytique de diff\'eomorphismes de $(\RR^{q_2},0)$. Une {\bf transformation semi-diagonale de param\`etres} $(t,\psi)$ ( et {\bf de puissances} $r_j$ et {\bf de taille} $k$), est un morphisme  $T_{t,\psi}: (x',\mu,\nu')\in ((\RR^{+*})^{k}\times\RR^{q_1}\times\RR^{q_2},0)\mapsto (x,\mu,\nu)\in  ((\RR^{+*})^{k}\times\RR^{q_1}\times\RR^{q_2},0)$ avec

$$
 x_1=t x'_1,\quad x_j=t^{r_1\times\cdots\times r_{j-1}}x'_j\quad\text{pour}\quad j=2,\ldots,k\quad\text{et}\quad \nu=\psi_\mu (\nu')
$$

\noindent Ces transformations forment un groupe pour la composition: $T_{t,\psi}\circ T_{t',\psi'}=T_{tt',\psi\circ\psi'}$. Ce groupe agit sur l'espace des couples d'Hilbert \`a composantes dans l'alg\`ebre $QR{\Cal H}^{k,q}_{cvg}$ (et même dans l'alg\`ebre $QR{\Cal H}^{k,q}$), de puissances $r_j$ et de dimension de non-trivialité $k-1$.\par

\vskip 3mm

\proclaim{Lemme 2.1.3} La multiplicit\'e alg\'ebrique du couple $(\chi,f)$ sur $B'_{k,0}$ est invariante par le groupe des transformations semi-diagonales.
\endproclaim

\vskip 3mm

\noindent{\bf Preuve.} Appliquons le premier \'eclatement au couple $(\chi,f)$ et \`a son image $(\chi',f')$ par $T_{t,\psi}$: $t^s\chi'=(T_{t,\psi}^{-1})_*\chi$ pour une certaine puissance $s$, et $f'=T_{t,\psi}^*(f)$. Soit $\widetilde{T}_{t,\psi}$ le relev\'e de $T_{t,\psi}$ au dessus des singularit\'es principales sur ${\Cal D}$ et ${\Cal D}'$, et soient $(\rho,\mu,\nu,u)$ et $(\rho',\mu,\nu',u')$ les coordonn\'ees au dessus de ${\Cal D}$ et ${\Cal D}'$. En utilisant les intégrales premières non triviales de $\chi$ et $\chi'$, dans ce premier éclatement, on vérifie facilement que

$$
 \rho=t\rho'\quad \nu=\psi_\mu(\nu')\quad\text{et}\quad u=u'
$$

\noindent Le morphisme $\widetilde{T}_{t,\psi}$ induit donc un isomorphisme $\widetilde{T}^*_{t,\psi}$ entre les anneaux $QR{\Cal H}^{1,q}(\rho,$\par
\noindent $\alpha,u)$ et $QR{\Cal H}^{1,q}(\rho',\alpha',u')$ qui pr\'eserve les degr\'es des blocs $\cdot$-homog\`enes et donc la multiplicit\'e alg\'ebrique.\qed

\vskip 3mm

\noindent{\bf Remarque 2.1.1.}  Soit maintenant $B_{k,0}=\{\prod_{j=1}^k x_j=0,\ \mu=0\}\cap\overline{U}$. Soit $m\in B_{k,0}\setminus B'_{k,0}$, d'après le procédé de désingularisation de la dérivation d'Hilbert (IVA[Mo]), le germe en $m$ du couple $(\chi,f)$ est {\bf analytiquement conjugu\'e \`a un couple d'Hilbert} qu'on note $(\chi_m,f_m)$, et qui est  \`a composantes dans une alg\`ebre $QR{\Cal H}_{cvg}^{k',q'}$ avec $k'<k$. On d\'efinit alors la multiplicit\'e alg\'ebrique du couple $(\chi,f)$ en $m$ par $ma(\chi,f)_m=ma(\chi_m,f_m)_m$. De plus, d'après ce même procédé, si $m\neq 0$ et $b=\pi^{-1}(m)$, alors les couples $(\chi_m,f_m)$ et $(\widetilde{\chi}_b,\widetilde{f}_b)$ sont conjugu\'ees par une transformation semi-diagonale $T_{\rho',\psi}$. La g\'en\'eralisation du lemme 2.1.1 est\par

\vskip 3mm

\proclaim{Lemme 2.1.4} L'application $m\in B_{k,0}\mapsto ma(\chi,f)_m$ est localement major\'ee.
\endproclaim

\vskip 3mm

\noindent{\bf Preuve.} Par induction sur $k$. Le lemme 2.1.1 implique le cas $k=1$. Soit $k>1$, prouvons le lemme au point $m=0$ par exemple. Soit $a\in \partial{\Cal D}$, le couple relev\'e $(\widetilde{\chi}_a,\widetilde{f}_a)$ est un couple d'Hilbert, \`a composantes dans une alg\`ebre $QR{\Cal H}^{k'+1,q'}_{cvg}$ avec $1\leq k'<k$ ($k'-1$ étant sa dimension de non-trivialité). Soit $B_{k',0}(a)\subset\overline{{\Cal N}}$ le germe en $a$ de l'intersection du bord $\partial {\Cal N}$ avec le sous-ensemble $\{\mu=0\}$. Soit $b\in B_{k',0}(a)\setminus(\partial{\Cal D},a)$ et soit $m=\pi(b)\in B_{k,0}\setminus B'{k,0}$. D'apr\`es la remarque 2.1.1, le germe en $b$ du couple $(\widetilde{\chi}_a,\widetilde{f}_a)$ est analytiquement conjugué au couple d'Hilbert $(\widetilde{\chi}_b,\widetilde{f}_b)$. De plus, les couples $(\chi_m,f_m)$ et $(\widetilde{\chi}_b,\widetilde{f}_b)$ sont conjugu\'ees par une transformation semi-diagonale $T_{\rho',\psi}$. Par le lemme 2.1.3, elles ont la m\^eme multiplicit\'e algébrique en $m$ et en $b$. Par l'hypoth\`ese de r\'ecurrence, cette multiplicit\'e algébrique est localement major\'ee en $a$, et par la compacit\'e de $\partial{\Cal D}$ et la surjectivité de $\pi$, elle est major\'ee sur un voisinage de 0 dans $B_{k,0}\setminus B'_{k,0}$. Le lemme 2.1.2 permet de finir la preuve.\qed

\vskip 3mm

 Je ne sais pas si cette application est semi-continue sup\'erieurement. Pour cela, il faut savoir comparer la multiplicit\'e algébrique principale $ma(\chi,f)_0$ et les multiplicit\'es algébriques secondaires $ma(\widetilde{\chi}_a,\widetilde{f}_a)_a$ sur le bord $\partial{\Cal D}$. Une généralisation de la notion de multiplicité algébrique, qui convient à notre contexte, est la suivante: soit $m\in B'_{k,0}$, on définit, par récurrence sur $k$, la {\bf multiplicité algébrique généralisée} du couple $(\chi,f)$ en $m$, qu'on note $mag(\chi,f)_m$: pour $k=1$, on pose $mag(\chi,f)_m=ma(\chi,f)_m$, et pour $k>1$, on pose

$$
mag(\chi,f)_{m}=\max\{\sup\{mag(\widetilde{\chi}_{a},\widetilde{f}_{a})_{a};\ a\in\partial {\Cal D}(0,\nu_m) \},\ ma(\chi,f)_m\}
$$

\noindent Pour $m\in B_{k,0}\setminus B'_{k,0}$, on pose $mag(\chi,f)_m=mag(\chi_m,f_m)_m$.\par

\vskip 3mm

\proclaim{Lemme 2.1.5} La multiplicité algébrique g\'en\'eralis\'ee est finie et est invariante dans les transformations semi-diagonales. L'application $m\in B_{k,0}\mapsto mag(\chi,f)_{m}$ est semi-continue sup\'erieurement.
\endproclaim

\vskip 3mm

\noindent{\bf Preuve.} Par une récurrence sur $k$. Le cas $k=1$ découle du lemme 2.1.1 (pour la semi-continuité supérieure et la finitude), et du lemme 2.1.3 (pour l'invariance par les transformations semi-diagonales) . Soit $k>1$, faisons la preuve en $m=0$ par exemple. Reprenons la transformation $T_{t,\psi}$ du lemme 2.1.3, en utilisant l'expression du morphisme $\pi$, on v\'erifie ais\'ement que sur le bord de ${\Cal D}$, cette transformation se rel\`eve en une transformation semi-diagonale $T_{1,\psi_t}$ avec $\psi_{t,\mu'}(\nu',X(\rho'_0,\mu'))$ $=(\nu,X(\rho_0,\mu'))$, o\`u $X(.)$ sont les fonctions \'el\'ementaires des alg\`ebres locales $QR{\Cal H}^{.,.}$ correspondantes. Donc, par l'hypothèse de récurrence, la multiplicité algébrique généralisée est invariante par les transformations semi-diagonales.\par

\vskip 3mm

  Reprenons les notations de la preuve du lemme 2.1.4. Par l'hypothèse de récurrence, la multiplicité algébrique généralisée est finie sur $B_{k',0}(a)$, et y est semi-continue supérieurement. Par la compacité de $\partial {\Cal D}(0)$, elle est finie et semi-continue supérieurement sur $\partial{\Cal N}\cap\{\mu=0\}$ (germifié le long de $\partial{\Cal D}(0)$). Donc, par son invariance par les transformations semi-diagonales, elle est finie et semi-continue supérieurement en $m=0$.\par\qed

\vskip 3mm

 La g\'en\'eralisation de cette notion de multiplicit\'e algébrique aux points $m$ tels que $\mu_m\neq 0$ est d'int\'er\^et th\'eorique certain mais n'est d'aucune utilit\'e dans notre contexte. En fait, cette multiplicit\'e algébrique est li\'ee aux blocs $\widetilde{\chi}_m$-homog\`enes de la s\'erie de $\widetilde{f}_m$, et ces blocs d\'ependent \'etroitement des valeurs $\mu_m$ par l'interm\'ediaire des valeurs propres des op\'erateurs d'Euler associ\'es \`a la dérivation  $\widetilde{\chi}_m$.\par

\vskip 3mm

\noindent{\bf (b) Couples d'Hilbert} $r${\bf -alg\'ebrisables et voisinages universels.}\par

\vskip 3mm 

 Soit $(\chi,f)$ un couple d'Hilbert \`a composantes $D_1,\ldots,D_{k-1},f$ dans l'alg\`ebre $QR{\Cal H}^{k+p,q}(x,\rho,\alpha)$. Soient $g_j=x_j^{r_j}(1+D_j)-x_{j+1}$, $j=1,\ldots,k-1$, les int\'egrales premi\`eres non triviales de $\chi$. On appelle {\bf partie principale} de $\chi$, qu'on note $\chi_{pr}$, la d\'erivation d'Hilbert dont les int\'egrales premi\`eres non triviales sont $g_{j,pr}=x_j^{r_j}-x_{j+1}$, $j=1,\ldots,k-1$.\par

\vskip 3mm

 Soient $X_j(x_j,\mu)=(x_j,z(x_j,\mu))$ et $X(x,\mu)=(X_j(x_j,\mu))_{j=1,\ldots,k}$ les fonctions \'el\'ementaires de l'alg\`ebre $QR{\Cal H}^{k+p,q}$ dans les coordonn\'ees $x$. Le couple $(\chi,f)$ est dit 0-{\bf alg\'ebrisable de degr\'e} $m_0\in\NN$ si les fonctions $D_j$ (resp. $f$) sont alg\'ebriques en $X_j$ (resp. $X$), de degr\'e $m_0$.\par

\vskip 3mm

 Soit $(\pi,{\Cal N})$ le premier \'eclatement de $\chi$ en 0 de diviseur exceptionnel $\overline{{\Cal D}}$. Soit $a\in\partial {\Cal D}$ de codimension $k'<k$ dans $\overline{{\Cal D}}$, et soient $D'_1,\ldots,D'_{k'-1},f'$ les composantes du couple d'Hilbert relev\'e  $(\widetilde{\chi}_a,\widetilde{f}_a)$. Soient $(x'=(x'_1,\ldots,x'_{k'}),\rho',\alpha')$ les coordonn\'ees locales sur le germe $({\Cal N},a)$ et soient $X'_j(x'_j,\mu')$ et $X'(x',\mu')$ les fonctions \'el\'ementaires, dans les coordonn\'ees $x'$, de l'alg\`ebre locale en $a$ $QR{\Cal H}^{k'+p+1,q'}$. Soit $m\in\NN$, on d\'efinit une application $EP_{0,a,m}$ ({\bf E}clatement puis {\bf P}erturbation) de l'ensemble des couples d'Hilbert de l'alg\`ebre $QR{\Cal H}^{k+p,q}$ dans l'ensemble des couples d'Hilbert de l'alg\`ebre $QR{\Cal H}^{k'+p+1,q'_m}$, de la fa\c con suivante: soient $D"_1,\ldots,D"_{k'-1}$ (resp. $f"$) des polyn\^omes en $X'_j$ (resp. $X'$), de degr\'e $m$, et \`a {\bf coefficients universels dans un voisinage de 0}. On note

$$
EP_{0,a,m}(\chi,f)_0=({\Cal X},F)_a
$$

\noindent le couple d'Hilbert de composantes $D'_1+D"_1,\ldots,D'_{k'-1}+D"_{k'-1}$ et $f'+f"$; les d\'erivations $\widetilde{\chi}_a$ et ${\Cal X}$ ayant les m\^emes parties principales, et $q'_m=q'+(0,q'_{2,m})$, où $q'_{2,m}$ est le nombre des coefficients universels dans les polynômes $D"_j$ et $f"$. Si $k'=1$, les fonctions $D'_j$ et $D"_j$ sont nulles. Ces applications envoient bien s\^ur une alg\`ebre convergente sur une alg\`ebre convergente.\par

\vskip 3mm

 Soient $r,m_0,\ldots,m_r,m_{r+1}\in\NN$. Supposons d\'efinis les couples d'Hilbert $r${\bf -alg\'ebrisables de degr\'e} $(m_0,\ldots,m_r)$. Un couple d'Hilbert $(\chi,f)_a$ est dit $(r+1)${\bf -alg\'ebrisable de degr\'e} $(m_0,\ldots,m_r,m_{r+1})$ s'il existe un couple $(\chi',f')_{a'}$ $r$-alg\'ebrisa\-ble de degr\'e $(m_0,\ldots,m_r)$ telle que

$$
(\chi,f)_a=EP_{a',a,m_{r+1}}(\chi',f')_{a'}
$$

\noindent Autrement dit, on obtient le couple $(\chi,f)_a$ en éclatant le couple $(\chi',f')$ en $a'$, puis en perturbant le germe en $a$ de son éclaté, par des fonctions "algébriques" de degré $m_{r+1}$.\par

\vskip 3mm

\noindent{\bf (b1) Construction de la collection universelle de couples} $r${\bf -alg\'ebrisables.}\par

\vskip 3mm

 Soit $(\chi_0,f_0)$ un couple d'Hilbert dont les composantes $D_1,\ldots,D_{k-1},f_0$ sont des \'el\'ements de $QR{\Cal H}^{k,q=(q_1,q_2)}$. D'apr\`es la d\'esingularisation de la d\'erivation d'Hilbert (cf. IVA, [Mo]), le diviseur exceptionnel de la d\'esingularisation compl\`ete de $\chi_0$ ($k-1$ \'eclatements) est identique \`a celui de la d\'esingularisation compl\`ete de la d\'erivation principale associ\'ee. Les strates de ce diviseur total (qui est compact et qu'on notera ${\Cal D}$ dans la suite), sont des {\bf semi-alg\'ebriques} dont le nombre et la complexit\'e ne d\'ependent que de l'entier $k$. Par exemple, la première de ces strates ${\Cal D}_{1,k}(0)$ (celle de plus grande dimension qui est $k-1$), est l'image de la semi-sphère
 
$$
S^{+*}_k(0,1)=\{y\in(\RR^{+*})^k;\ y_1+\cdots y_k=1\}
\tag 9
$$

\noindent par le morphisme ${\Cal L}_{k,0}: y\in S^{+*}_k(0,1)\mapsto u\in\RR^{k-1}$, avec $u_j=y_j-y_{j+1}$. Les autres strates ${\Cal D}_{i,k}(0)$ (pour $i>1$), sont chacune identique à une strate ${\Cal D}_{1,k'}(0)$ ($k'<k$), ou à un produit de tels strates. Pour simplifier la présentation dans la suite, on omettra les indices $k$ et $(0)$, sauf en cas d'ambiguité.\par

\vskip 3mm

 On suppose que le couple $(\chi_0,f_0)$ est 0-alg\'ebrisable de degr\'e $m_0$, \`a {\bf coefficients universels dans un voisinage de 0}. Les {\bf propri\'et\'es de finitude} de ce couple ($\chi_0$-r\'egularit\'e et $\chi_0$-finitude de $f_0$ (cf. IB, [Mo])), ne d\'ependent que des entiers $k$ (nombre de variables $x$), $q_1$ (nombre de fonctions élémentaires dans les variables $x$) et $m_0$ (degré dans ces variables et dans ces fonctions): la référence la plus générale sur ce sujet (et sur ce qui va suivre), est le travail de Khovanski [K], notament sur les fonctions pfaffiennens définies sur des variétés semi-pfaffiennes. Cependant, certains travaux récents d'Il'yashenko, Novikov et Yakovenko ([IY], [NY]), consacrés à l'étude de l'intégrale Abélienne, pourraient aider à donner des estimations explicites des degrés de régularité et des indices de noethérianité.\par
 
\vskip 3mm

  A partir de ce couple, nous allons construire sur ${\Cal D}$ une {\bf collection universelle de couples} $r${\bf -alg\'ebrisables} et une {\bf collection associ\'ee de voisinages universels} (en ce sens qu'ils ne d\'ependent que des entiers $m_0,k$ et $q_1$), et au dessus de laquelle on a des {\bf propri\'et\'es de finitudes universelles} pour les couples d'Hilbert $(\chi,f)$ "voisins" de $(\chi_0,f_0)$ (dans un sens qu'on pr\'ecisera).\par

\vskip 3mm

 Soit $m_1$ un majorant de la multiplicit\'e alg\'ebrique g\'en\'eralis\'ee   $mag(\chi_0,f_0)_0$. Soit $(\pi_1,{\Cal N}_1)$ le premier \'eclatement de $\chi_0$ en 0, de diviseur exceptionnel $\overline{{\Cal D}}_1$. Soit $(\chi_1,f_1)$ la collection de  couples d'Hilbert au dessus de $\partial{\Cal D}_1$, qui sont 1-alg\'ebrisables de degr\'e $(m_0,m_1)$, et qui sont construits comme suit: en tout point $a\in\partial{\Cal D}_1$, on pose

$$
(\chi_1,f_1)_a=EP_{0,a,m_1}(\chi_0,f_0)_0
$$

\noindent Notons $(\widetilde{\chi}_0,\widetilde{f}_0)$ le relev\'e du couple $(\chi_0,f_0)$ et $(\widetilde{\chi}_0,\widetilde{f}_0)_a$ son germe en $a\in\overline{{\Cal D}}_1$. Soit ${\Cal C}_{1,\ell,i}$ une strate de $\partial{\Cal D}_1$ de codimension $\ell\in\{1,\ldots,k-1\}$ dans $\overline{{\Cal D}}_1$ (l'indice $i$ est \'enum\'eratif, et il est majoré par une fonction de $k$)). C'est un semi-algébrique de type (9) (en remplaçons $k$ par $k-\ell$). Le long de cette strate, la dimension de non trivialit\'e de la d\'erivation $\widetilde{\chi}_0$ est $\ell-1$. Pour pouvoir construire des couples 2-alg\'ebrisables {\bf universels}, on doit montrer que la multiplicit\'e alg\'ebrique g\'en\'eralis\'ee des couples 1-alg\'ebrisables construits, est major\'ee sur $\partial{\Cal D}_1$. Ceci est {\bf hautement non trivial}, car ces couples ne se recollent pas forc\'ement d'une strate \`a l'autre, et les seules strates fermées sont les strates réduites à un point.\par

\vskip 3mm

  Soit $a\in\partial{\Cal C}_{1,\ell,i}\cap {\Cal C}_{1,\ell_1,i_1}$ (avec forcément $\ell_1>\ell$). Notons simplement ${\Cal C}_a$ le germe de ${\Cal C}_{1,\ell,i}$ en $a$. Si $b\in{\Cal C}_a$, on sait que (cf. remarque 2.1.1), le germe en $b$ du couple d'Hilbert $(\widetilde{\chi}_0,\widetilde{f}_0)_a$ est analytiquement conjugu\'e au couple d'Hilbert $(\widetilde{\chi}_0,\widetilde{f}_0)_b$. Dans la suite, on parlera indiff\'erement de l'un ou de l'autre dans les situations similaires. Appliquons un premier \'eclatement $(\pi_{1,a},{\Cal N}_{1,a})$ au couple $(\widetilde{\chi}_0,\widetilde{f}_0)_a$, au point $a$, et notons $(\widetilde{\chi}_{0,1},\widetilde{f}_{0,1})$ son relev\'e. Soit $a_1\in \pi_{1,a}^{-1}(a)\cap \overline{\pi_{1,a}^{-1}({\Cal C}_a)}$, et soit ${\Cal C}_{a_1,a}$ le germe en $a_1$ de $\pi_{1,a}^{-1}({\Cal C}_a)$. Par la remarque 2.1.1, en tout point $b_1\in {\Cal C}_{a_1,a}$, le couple $(\widetilde{\chi}_0,\widetilde{f}_0)_{\pi_{1,a}(b_1)}$ est \'equivalent au germe en $b_1$ du couple $(\widetilde{\chi}_{0,1},\widetilde{f}_{0,1})_{a_1}$ par une transformation semi-diagonale $T_{t_1,\psi_1}$.\par

\vskip 3mm

 On r\'ep\`ete ce proc\'ed\'e en $a_0=a$ au plus $p=\ell_1-\ell$ fois, appliqu\'e aux singularit\'es $a_j$ de dimension de non trivialit\'e$>\ell$. On construit ainsi un morphisme $(\pi_{p,a},{\Cal N}_{p,a})$ et un relev\'e $(\widetilde{\chi}_{0,p},\widetilde{f}_{0,p})$ du couple 
$(\widetilde{\chi}_0,\widetilde{f}_0)_a$. Soit $a_p\in \pi_{p,a}^{-1}(a)\cap \overline{\pi_{p,a}^{-1}({\Cal C}_a)}$ et soit ${\Cal C}_{a_p,a}$ le germe en $a_p$ de $\pi_{p,a}^{-1}({\Cal C}_a)$. Soit $b_p\in {\Cal C}_{a_p,a}$ et $b=\pi_{p,a}(b_p)$. Le couple $(\widetilde{\chi}_0,\widetilde{f}_0)_b$ est \'equivalent au germe en $b_p$ du couple $(\widetilde{\chi}_{0,p},\widetilde{f}_{0,p})_{a_p}$ par une transformation semi-diagonale $T_{t_p,\psi_p}$. De plus, le germe en tout point du couple $(\widetilde{\chi}_{0,p},\widetilde{f}_{0,p})$ est de dimension de non trivialit\'e $\ell-1$. Construisons alors en $a_p$ un couple $(\chi_{1,\ell,i},f_{1,\ell,i})_{a_p}$ 
au plus $(p+1)$-alg\'ebrisable de degr\'e $(m_0,0,\ldots,0,m_1)$, \`a partir du couple $(\chi_{0},f_{0})$ et de polyn\^omes universels perturbateurs $D_{j,\ell,i}$ et $F_{\ell,i}$, \'ecrits dans les coordonn\'ees locales en $a_p$. {\bf Le point cl\'e} est que le couple $(\chi_1,f_1)_b$ est \'equivalent au germe en $b_p$ du couple 
$(\chi_{1,\ell,i},f_{1,\ell,i})_{a_p}$ par une transformation semi-diagonale $T_{t_p,\psi'_p}$ o\`u $\psi'_p$ est le composé de $\psi_p$ et d'un automorphisme de l'espace des coefficients universels des polyn\^omes $D_{j,\ell,i}$ et $F_{\ell,i}$.\par

\vskip 3mm

 En appliquant cette d\'esingularisation en tout point $a$ du bord de la strate ${\Cal C}_{1,\ell,i}$, on construit une vari\'et\'e compact, connexe \`a bord $\overline{{\Cal D}}_{1,\ell,i}$ (qui est une union de strates de type (9), dont la complexité est majorée par l'entier $k$, et dont la strate principale ${\Cal D}_{1,\ell,i}$ est isomorphe \`a ${\Cal C}_{1,\ell,i}$), sur laquelle est d\'efini un couple local $(\chi_{1,\ell,i},f_{1,\ell,i})$, dont la classe d'\'equivalence (ou classe de recollement) est le groupe des transformations semi-diagonales. Par le lemme 2.1.5, la multiplicit\'e alg\'ebrique g\'en\'eralis\'ee des couples $(\chi_1,f_1)$, est donc major\'ee sur la strate ${\Cal C}_{1,\ell,i}$, et par conséquent, elle est majorée sur $\partial{\Cal D}_1$, par une fonction universelle dans les entiers $m_0,m_1,k$ et $q_1$.\par

\vskip 3mm

 Soit $m_2$ un majorant de cette multiplicit\'e algébrique généralisée sur $\partial{\Cal D}_1$. Appliquons un premier \'eclatement (dont le germe en $a$ est) $(\pi_2,{\Cal N}_2)_a$ aux couples $(\chi_1,f_1)_a$ en tout point $a\in\partial{\Cal D}_1\setminus\cup{\Cal C}_{1,1,i}$. Soit $\overline{{\Cal D}}_2$ son diviseur exceptionnel (il n'est pas forcément connexe). On construit alors une collection de couples d'Hilbert $(\chi_2,f_2)$ au dessus de $\partial{\Cal D}_2$, qui sont 2-alg\'ebrisables de degr\'e $(m_0,m_1,m_2)$: en tout point $a_1\in\partial{\Cal D}_2\cap\pi_2^{-1}(a)$, on pose

$$
(\chi_2,f_2)_{a_1}=EP_{a,a_1,m_2}(\chi_1,f_1)_a
$$

\vskip 3mm

 On r\'ep\`ete ce proc\'ed\'e au plus $(k-1)$-fois appliqu\'e aux singularit\'es de dimension de non trivialit\'e$>0$. Pour majorer la multiplicit\'e alg\'ebrique g\'en\'eralis\'ee \`a une \'etape d'ordre $r$, on plonge chaque strate ${\Cal C}_{r,\ell,i}$ dans la vari\'et\'e compacte correspondante $\overline{{\Cal D}}_{r,\ell,i}$ sur laquelle on construit un couple local $(\chi_{r,\ell,i},f_{r,\ell,i})$ \`a partir des couples locaux $(\chi_{r',\ell',i'},f_{r',\ell',i'})$ sur les vari\'et\'es compactes ant\'erieures $\overline{{\Cal D}}_{r',\ell',i'}$.\par

\vskip 3mm

\noindent{\bf (b2) Construction de la collection des voisinages universels.}\par

\vskip 3mm

 Soit $(\chi_0,f_0)$ comme ci-dessus et $m_1,\ldots,m_{k-1}$ des majorants des multiplicit\'es alg\'ebriques g\'en\'eralis\'ees successives. Soit $((\chi_r,f_r))_{r=1,\ldots,k-1}$ la collection des couples $r$-alg\'ebrisables associ\'ee. Les diff\'erentes vari\'et\'es ${\Cal D}$, ${\Cal D}_r$ et ${\Cal D}_{r,\ell,i}$ sont des ensembles semi-alg\'ebriques dont la complexit\'e ne d\'epend que de l'entier $k$. En vue d'obtenir des propriétés de finitude universelles pour la collection $((\chi_r,f_r))_{r=0,\ldots,k-1}$, on construira, dans la suite, des voisinages {\bf au-dessus} de ces vari\'et\'es, en consid\'erant le param\`etre $\nu$ dans un voisinage de 0 dans l'espace de la collection des coefficients universels des couples $((\chi_r,f_r))_{r=0,\ldots,k-1}$. La dimension de cet espace est {\bf universelle} dans les entiers $(m_j),k$ et $q_1$.\par

\vskip 3mm

 Pla\c cons-nous au dessus d'une vari\'et\'e $\overline{{\Cal D}}_{r,\ell,i}$. Soit $V_{r,\ell,i}(a)\subset\overline{{\Cal D}}_{r,\ell,i}$ un voisinage {\bf choisi} pour la multiplicit\'e alg\'ebrique g\'en\'eralis\'ee du couple $(\chi_{r,\ell,i},f_{r,\ell,i})_a$: autrement dit, cette multiplicit\'e en tout point de $V_{r,\ell,i}(a)$ est major\'ee par la multiplicité en $a$ (cf. lemme 2.1.5). Les anneaux $QR{\Cal H}_{cvg}$ sont noeth\'eriens. D'apr\`es le lemme de coh\'erence IB3[Mo], l'indice de noeth\'erianit\'e des id\'eaux de faisceaux diff\'erentiels d'un anneau noeth\'erien, est une fonction semi-continue sup\'erieurement. {\bf Choisissons} un voisinage $W_{r,\ell,i}(a)$ {\bf au dessus} de $V_{r,\ell,i}(a)$, dans lequel cet indice est major\'e par celui de la fibre en $a$, et qu'on note $L_{r,\ell,i}(a)$. Soit $L_{r,\ell,i}((m_j),k,q_1)$ son maximum sur $\overline{{\Cal D}}_{r,\ell,i}$ (il est invariant par les transformations semi-diagonales!). Et soit $L((m_j),k,q_1)$ son maximum sur la collection des couples $cdot$-alg\'ebrisables.\par

\vskip 3mm

 Appliquons un premier \'eclatement $(\pi_a,{\Cal N}_a)$ en $a$ de diviseur exceptionnel $\overline{{\Cal D}}_{r+1,a}$; sa strate principale ${\Cal D}_{r+1,a}$ est incluse dans une composante connexe de ${\Cal D}_{r+1}$, et elle ne d\'epend pas du point $a$ (en tant que fibre d'une fbration triviale). Notons simplement $(\widetilde{\chi},\widetilde{f})$ le relev\'e du couple $(\chi_{r,\ell,i},f_{r,\ell,i})_a$ au dessus de $\overline{{\Cal D}}_{r+1,a}$. Soient $j_1,j_2\in\{0,\ldots,L_{r,\ell,i}(a)\}$ et soient

$$
f_{j_1,j_2}^{\pm}=(\widetilde{\chi}^{j_1}\pm\widetilde{\chi}^{j_2})\widetilde{f}
\tag 10
$$

\noindent Ce sont les fonctions (ou leurs germes) qui apparaissent dans la preuve du lemme de finitude IB1[Mo]. Et il s'agit donc de montrer que le degr\'e de la projection intégrale $\pi_{\widetilde{\chi}}$ restreinte \`a la fibre nulle de $f_{j_1,j_2}^{\pm}$, est major\'e au dessus de ${\Cal D}_{r+1,a}$. D'apr\`es le lemme d'extension IB2[Mo], l'alg\`ebre $QR{\Cal H}^{1,.}_{cvg}$ est $\widetilde{\chi}$-finie. Ce degr\'e est donc major\'e au dessus de tout compact de ${\Cal D}_{r+1,a}$.\par

\vskip 3mm

  Malheureusement, ces fonctions $f_{j_1,j_2}^{\pm}$ ne se prolongent pas forc\'ement au bord de ${\Cal D}_{r+1,a}$. Pour rem\'edier \`a cela, il vaut mieux remplacer la d\'erivation $\widetilde{\chi}$ par la d\'erivation $\chi_{r,\ell,i}$ (ce qui est tr\`es convenant pour les estimations explicites), ou simplement par la d\'erivation ${\Cal Y}$ donn\'ee par la proposition IVA1[Mo]

$$
({\Cal L}_\ell)_*{\Cal Y}=F_\ell\widetilde{\chi}
$$

\noindent le morphisme ${\Cal L}_\ell$ et la fonction $F_\ell$ ne d\'ependent que de la d\'erivation $\chi_{r,\ell,i}$. De plus, la d\'erivation ${\Cal Y}$ et la fonction $F_\ell$ admettent un prolongement au bord de la semi-sph\`ere $S^{+*}_\ell(0,1)=({\Cal L}_{\ell,0})^{-1}({\Cal D}_{r+1,a})$ (qui est du type (9)). Maintenant, les fonctions

$$
F_{j_1,j_2}^{\pm}=F_\ell^{2L_{r,\ell,i}(a)}(f_{j_1,j_2}^{\pm}\circ{\Cal L}_\ell)
$$

\noindent admettent aussi un prolongement au bord de $S^{+*}_\ell(0,1)$, et leur degr\'e dans la projection $\pi_{\Cal Y}$ majore le degr\'e recherch\'e. Toujours d'apr\`es le lemme d'extension IB2[Mo], les alg\`ebres $QR{\Cal H}_{cvg}$ sont ${\Cal Y}_{A}$-finie pour tout $A\in \overline{S^{+*}_\ell(0,1)}$. Notons $d_{r,\ell,i}(a)$ la somme de ces degr\'es pour les fonctions $f_{j_1,j_2}^{\pm}$, au dessus de ${\Cal D}_{r+1,a}$. C'est une fonction semi-continue sup\'erieurement sur $\overline{{\Cal D}}_{r,\ell,i}$. {\bf Choisissons} un voisinage $U_{r,\ell,i}(a)\subset W_{r,\ell,i}(a)$ tel que ce nombre soit r\'ealis\'e sur $\pi_a^{-1}(U_{r,\ell,i}(a))$. Notons $d_{r,\ell,i}((m_j),k,q_1)$ le maximum de $d_{r,\ell,i}(a)$ sur $\overline{{\Cal D}}_{r,\ell,i}$. Ainsi, la collection des voisinages $U_{r,\ell,i}(a)$ est {\bf universelle} dans les entiers $(m_j),k$ et $q_1$.\par

\vskip 3mm

\noindent{\bf (c) Preuve du théorème 2.1.}\par

\vskip 3mm

 Pour éviter un conflit de notations dans la suite, on notera $({\Cal X},F)$ au lieu du couple d'Hilbert $(\chi,f)$ du théorème 2.1. Pour alléger le texte dans la suite, on entendra par {\bf restriction appropri\'ee}, toute restriction \`a un sous-ensemble de $W_0$ ou $U_0$ ou \`a leurs relev\'es dans la d\'esingularisation de ${\Cal X}$. La preuve qui suit est bas\'ee sur les r\'esultats de la partie IVB de [Mo]. Le couple $({\Cal X},F)$ satisfait à l'hypothèse $(H\lambda)$. Par le théorème IVB1[Mo], il possède une multiplicité algébrique $ma=ma({\Cal X},F)_0$. Prenons $m_0=ma+N_0+2$. Soit $(\chi_0,f_0)$ le couple 0-alg\'ebrisable de degr\'e $m_0$, \`a composantes dans $QR{\Cal H}^{k,(q_1,.)}(x,\alpha')$ ($\alpha'=(\mu,\nu')$, cf. ci-desous pour les variables $\nu'$). Construisons une collection universelle $((\chi_r,f_r))_{r=1,\ldots,k-1}$ de couples $r$-alg\'ebrisables de degr\'e $(m_0,\ldots,m_r)$ avec 

$$
m_r=\sup_{a\in\partial{\Cal D}_{r-1}}\{mag(\chi_{r-1},f_{r-1})_a\}+2
$$

\noindent Ainsi, les degr\'es $m_r$ ne d\'ependent que des entiers  $m_0,k$ et $q_1$. Les variables $\nu'$ étant donc dans l'espace de la collection des coefficients universels de ces couples $r$-algébrisables. Notons 

$$
L'_{r,\ell,i}(m_0,k,q_1)=L_{r,\ell,i}((m_j),k,q_1),\quad L'(m_0,k,q_1)=L((m_j),k,q_1)
$$

$$
d'_{r,\ell,i}(m_0,k,q_1)=d_{r,\ell,i}((m_j),k,q_1)
$$

\vskip 3mm

\noindent et notons toujours $U_{r,\ell,i}(a)$ les voisinages universels associés à ces degrés. Montrons que les propri\'et\'es de finitude du couple $({\Cal X},F)$ ne d\'ependent que de celles de cette collection. Appliquons une transformation semi-diagonale $T_{t,id}$ au couple $({\Cal X},F)$ et notons $({\Cal X}_t,F_t)$ son image. Soient $\lambda'=(\lambda'_1,\ldots,\lambda'_{k-1})$ les valeurs (coordonnées) des intégrales premières non triviales de ${\Cal X}_t$. Soit $\gamma$ l'orbite principale de ${\Cal X}$ en 0, son image $\gamma'=T^{-1}_{t,id}(\gamma)$ est l'orbite principale de ${\Cal X}_t$ en 0. Soit $\sigma$ (resp. $\sigma'=T^{-1}_{t,id}(\sigma)$) une transversale analytique à $\gamma$ (resp. $\gamma'$), de coordonnées $(\alpha,\lambda)$ (resp. $(\alpha,\lambda')$). La restriction de $T_{t,id}$ à $\sigma'$ s'écrit

$$
(\alpha,\lambda)=T_{t,id|\sigma'}(\alpha,\lambda')=(\alpha,(t^{r_1\times\cdots\times r_j}\lambda'_j))
$$

\noindent Donc, le morphisme $T_{t,id|\sigma'}^*$ envoie l'idéal ${\Cal M}_\lambda\subset\RR\{\alpha,\lambda\}$ sur l'idéal ${\Cal M}_{\lambda'}\subset\RR\{\alpha,\lambda'\}$. Par le lemme de transfert IB6[Mo], le couple $({\Cal X}_t,F_t)$ satisfait à l'hypothèse $(H\lambda)$, l'entier $N_0$ étant préservé. Toujours par ce même lemme, les indices de neothérianité de ces couples, sont préservés. Et par le lemme 2.1.3, les couples $({\Cal X},F)$ et $({\Cal X}_t,F_t)$ ont la m\^eme multiplicité algébrique en 0. Ils ont donc les m\^emes propri\'et\'es de finitude. Et pour $t$ suffisament petit, les coefficients des jets d'ordre $m_0$ en $X(x,\mu)$ des composantes de $({\Cal X}_t,F_t)$ sont r\'ealis\'es sur le projeté correspondant du voisinage universel $U_{0,k,1}(0)$ (après bien sûr identification des coordonnées $x$ et $x'$).\par

\vskip 3mm

 Ce qui suit utilise abondamment les arguments de la preuve du théorème IVB1 de [Mo]. Omettons provisoirement l'indice $t$. Soit $(\pi,{\Cal N})$ le premier \'eclatement du couple $({\Cal X},F)$ de diviseur exceptionnel $\overline{{\Cal D}}_1$ et soit $(\widetilde{{\Cal X}},\widetilde{F})$ son relev\'e. Soient $(\rho,\alpha,u)$ les coordonn\'ees sur ${\Cal N}$ au dessus de ${\Cal D}_1$ et soit $a_0$ la singularit\'e principale sur ${\Cal D}_1$. Notons $G_1={\bold j}^{N_0}_u(\widetilde{F}_{a_0})$ et $G_2=\sum_{n\leq ma} g_n(G_1)$ o\`u $g_n(G_1)$ est le bloc $\widetilde{{\Cal X}}_{a_0}$-homog\`ene de degr\'e $n$. Le germe $\widetilde{F}_{a_0}$ est $\widetilde{{\Cal X}}_{a_0}$-\'equivalent au germe

$$
 G(({\Cal X},F)(t))=G_2+{\bold j}^{ma+N_0+1}_{X(\rho,\mu)}(\widetilde{F}_{a_0}-G_1)
$$

\noindent qui est alg\'ebrique en $X$ de degr\'e$\leq m_0$. Soient $(\rho',\alpha',u')$ les coordonn\'ees de l'\'eclat\'e de $(\chi_0,f_0)$ et soit $a'_0$ sa singularit\'e principale. Soit $G(\chi_0,f_0)$ le germe construit de la m\^eme fa\c con que $G({\Cal X},F)(t))$, \`a partir du couple $(\chi_0,f_0)$. Ces deux fonctions coincident par une identification des coordonn\'ees $(\rho,u)$ et $(\rho',u')$, sur une restriction appropri\'ee. Quitte \`a r\'eduire le voisinage $U_{0,k,1}(0)$, on peut supposer que les invariants universels de $G(\chi_0,f_0)$ sont r\'ealis\'es sur le relev\'e en $a'_0$ de $U_{0,k,1}(0)$. Soit $L_0(N_0,ma,k,q_1)$ l'indice de noeth\'erianit\'e de l'id\'eal diff\'erentiel $I_{
\widetilde{\chi}_{0,a'_0}
,G(\chi_0,f_0)}$ et soit 

$$
 d_0(N_0,ma,k,q_1)=\sum_{j_1,j_2\leq L_0} d°\pi_{\widetilde{\chi}_{0,a'_0}|Z((\widetilde{\chi}_{0,a'_0}^{j_1}\pm\widetilde{\chi}_{0,a'_0}^{j_2})G(\chi_0,f_0))}
$$

\noindent Ainsi, sur une restriction appropri\'ee

$$
 ind(I_{\widetilde{{\Cal X}}_{a_0},\widetilde{F}_{a_0}})\leq L_0\quad\text{et}\quad d°\pi_{\widetilde{{\Cal X}}_{a_0}|Z(\widetilde{F}_{a_0})}\leq L_0 d_0
$$
 Soit $g_0$ la fonction obtenue \`a partir de $\widetilde{f}_0$ par identification des coordonn\'ees $(\rho,u)$ et $(\rho',u')$. Au dessus d'une restriction appropri\'ee sur ${\Cal D}_1\setminus\{a_0\}$, la fonction $\widetilde{F}$ est $\widetilde{{\Cal X}}$-\'equivalente \`a $g_0$. Ainsi, sur cette restriction

$$
 ind({\Cal I}_{\widetilde{{\Cal X}},\widetilde{F}})\leq L'(m_0,k,q_1)\quad\text{et}\quad d°\pi_{\widetilde{{\Cal X}}|Z(\widetilde{F})}\leq L' d'_{0,k,1}(m_0,k,q_1)
$$

 Sur le bord de ${\Cal D}_1$, on utilise l'argument de r\'ecurrence suivant: notons $({\Cal X}_r,F_r)$ le d\'esingularis\'e de $({\Cal X},F)(t)$ sur $\partial{\Cal D}_r$ \`a une \'etape d'ordre $r$. Soit $a$ un point d'une restriction appropri\'ee sur $\partial{\Cal D}_r$, de coordonn\'ees $(x,\rho,\alpha)$, et soit $a'$ le point correspondant dans l'\'eclatement du couple universel $(\chi_{r-1},f_{r-1})$, de coordonn\'ees $(x',\rho',\alpha')$. Quitte \`a r\'eduire $(t,\rho)$, on peut supposer que la restriction appropri\'ee en $a$ est r\'ealis\'ee sur le voisinage universel $U_{r,\ell,i}(a')$. Soit $({\Cal X}_r,F_r)_a$ le relev\'e de $({\Cal X},F)$ et soit $\gamma_r\subset {\Cal D}_r$ l'orbite principale en $a$. Supposons que, apr\`es identification des coordonn\'ees $(x,\rho)$ et $(x',\rho')$ (et sur une restriction appropri\'ee), on ait

\roster

\item"$(i_r)$" il existe $N_{0,r}=(n_{0,r,j})$ tel que pour tout $m\in\gamma_r$, ${\Cal I}_{({\Cal X}_r,F_r)_a}(m)\supset (\rho^{N_{0,r}})$.

\item"$(ii_r)$" il existe $N_{1,r}=(n_{1,r,j})$ avec $n_{1,r,j}>n_{0,r,j}$ tel que le couple $({\Cal X}_r,F_r)_a$ soit \'equivalent au couple $(\widetilde{\chi}_{r-1},\widetilde{f}_{r-1})_{a'}$ dans le quotient $SB/(\rho^{N_{1,r}})$ (autrement dit, les composantes de ces deux couples d'Hilbert, sont équivalentes dans ce quotient).

\endroster

\noindent Montrons qu'il en est de m\^eme \`a l'\'etape d'ordre $r+1$. D'apr\`es l'hypoth\`ese $(ii_r)$ et la construction des couples $r$-alg\'ebrisables, les couples $({\Cal X}_r,F_r)_a$ et $(\chi_r,f_r)_{a'}$ sont \'equivalents dans le quotient $SB/(\rho^{N_{1,r}}){\Cal M}_x^{m_{r}}$ sur une restriction appropri\'ee. Appliquons un premier \'eclatement aux couples $({\Cal X}_r,F_r)_a$ et $(\chi_r,f_r)_{a'}$ de diviseurs exceptionnels $\overline{{\Cal D}}_{r+1}(a)$ et $\overline{{\Cal D}}_{r+1}(a')$. Soient $(\rho_0,\rho,\alpha,u)$ et $(\rho'_0,\rho',\alpha',u')$ les coordonn\'ees des \'eclat\'es. Soit $g_r$ la fonction obtenue \`a partir de $\widetilde{f}_r$ par identification des coordonn\'ees $(\rho_0,\rho,u)$ et $(\rho'_0,\rho',u')$. Au dessus d'une restriction appropri\'ee sur ${\Cal D}_{r+1}(a)$, la fonction $\widetilde{F}_r$ est $\widetilde{{\Cal X}}_r$-\'equivalente \`a $g_r$. Ainsi, sur cette restriction

$$
 ind({\Cal I}_{\widetilde{{\Cal X}}_r,\widetilde{F}_r})\leq L'\quad\text{et}\quad d°\pi_{\widetilde{{\Cal X}}_r|Z(\widetilde{F}_r)}\leq L' d'_{r,\ell,i}(m_0,k,q_1)
$$

\noindent Et d'apr\`es le th\'eor\`eme IVB1[Mo], et la construction des couples $(r+1)$-alg\'ebrisables, on a les hypoth\`eses $(i_{r+1})$ et $(ii_{r+1})$ avec

$$
 N_{0,r+1}=(m_r-1,N_{0,r}),\ N_{1,r+1}=(m_r,N_{1,r}),\ N_{0,1}=m_0-1\ \text{et}\ N_{1,1}=m_0
$$

 Soit $(\pi,{\Cal N})$ la d\'esingularisation compl\`ete de $({\Cal X}_t,F_t)$, de diviseur exceptionnel ${\Cal D}$ et soit $\Delta=\overline{\pi^{-1}(U_0)}\cap{\Cal D}$. Pour finir la preuve du lemme, on choisit une valeur de $t>0$ sur un recouvrement fini au dessus de $\Delta$, qui soit r\'ealis\'e dans les voisinages universels $U_{r,\ell,i}$. On obtient alors
 
$$
 ind(I_{{\Cal X},F})\leq \max\{L_0,L'\}\quad\text{et}\quad d°\pi_{{\Cal X}|Z(F)}\leq L_0d_0+L' \sum_{r=0,\ldots,k-1}\max_{\ell,i}\{d'_{r,\ell,i}\}
$$
\qed

\vskip 3mm

\noindent{\bf 2.2. Preuve du Lemme 2.}\par

\vskip 3mm

 Soit $(\chi,f)$ le couple d'Hilbert associ\'e au d\'eploiement $\omega_\nu$ au voisinage du cycle singulier $\Gamma_k$. On a $k=q_1$ qui est aussi le nombre de singularit\'es sur $\Gamma_k$; par le th\'eor\`eme de Bezout, il est major\'e par une fonction du degr\'e $d$. La restriction $W_0$ est un graphe analytique, donné par

$$
W_0=\{(\mu(\nu),\nu,\lambda(\nu);\ \nu\in (\RR^{(d+1)(d+2)},0)\}
$$

\noindent Or, l'id\'eal $\chi$-transverse de $f$ le long de $\gamma$, et l'hypothèse $(H\lambda)$ sont invariants dans les changements de coordonn\'ees analytiques sur les transversales \`a $\Gamma_k$ (lemme de transfert IB6[Mo]). Choisissons donc $x_j=t$ pour $\epsilon=0$. Toujours d'après ce lemme de transfert, l'idéal retriction $J_{\chi,f,\gamma|W_0}\subset\RR\{\nu\}$ coincide avec l'idéal $\partial/\partial x_1$-transverse de $p_{1,.}-id$ (cf. (7)), le long de l'orbite $\{\nu=0\}$. En conjuguant (7) à (8) par le difféomorphisme $x_1=t+O(\epsilon)$ on obtient que cet idéal coincide avec l'idéal $\partial/\partial t$-transverse de $p_{2,.}-id$ le long de $\{\nu=0\}$. D'après l'hypothèse $(**)$ et le lemme de division II1[Mo], on a donc $J_{\chi,f,\gamma|W_0}=(\epsilon)$. Or $\lambda_{|\epsilon=0}\equiv 0$, par cons\'equent $N_0=1$.\par

\vskip 3mm

   Appliquons un premier \'eclatement $(\pi,{\Cal N})$ au couple $(\chi,f)$. Soit $a_0$ la singularit\'e principale de $\widetilde{\chi}=\rho\partial/\partial\rho-\sum s_ju_j\partial/\partial u_j$. Comme $r_{j|\epsilon=0}\equiv 1$, on a $\rho=kt$ pour $\epsilon=0$. Par cons\'equent, en comparant les séries asymptotiques de $\widetilde{f}_{a_0}$ et de $p_{2,.}-id$, on obtient que les multiplicit\'es algébriques $ma(\chi,f)_0$ et $ma(I(t),0)$, coincident.\qed 

\vskip 5mm

  Parallélement à l'opération $EP_{0,a,m}$, définissons l'opération $PE_{0,a,m}: (\chi,f)_0\mapsto ({\chi'},f')_a$, qui consiste en une perturbation de degré $m$ des composantes du couple $(\chi,f)$, suivie d'un éclatement. Notons ${\Cal C}_0$ la classe des couples 0-algébrisables, à composantes dans les algèbres $QR{\Cal H}^{k,(q_1,.)}(x,\mu,.)$, de dimension de non trivialité $k-1$, et de puissances $r_j$ fixées. Si le couple $(\chi,f)_0\in{\Cal C}_0$ est de degré $m_0$, on dit que le couple image $(\chi',f')_a=PE_{0,a,m}(\chi,f)_0$ est 1'-algébrisable de degré $(m_0,m)$. Notons ${\Cal C}_0(q_2)$ la sous-classe des couples 0-algébrisables {\bf induits}: leurs composantes appartiennent à l'algèbre $QR{\Cal H}^{k,(q_1,q_2)}(x,\mu,\nu)$. Soit ${\Cal C}_1(m)=EP_{0,a,m}({\Cal C}_0(q_2))$ (resp. ${\Cal C}'_1(m')=PE_{0,a,m'}({\Cal C}_0(q_2))$) la sous-classe des couples 1-algébrisables (resp. 1'-algébrisables) {\bf induits}: leurs composantes appartiennent à une algèbre $QR{\Cal H}^{k',(q'_1,q_{2,m})}(y,\mu',\nu')$ (resp. $QR{\Cal H}^{k',q'_1,q_{2,m'})}(y,\mu',\nu")$). En vue de faciliter la recherche d'une estimation explicite dans le théorème 2.1, il est utile de s'intéresser à la question suivante: existe-il un entier $m'=m'(k,q_1,q_2,m)$, et un morphisme analytique (induction) $\psi(y,\mu',\nu')=(y,\mu',\nu"(\nu'))$, tels que le diagramme suivant soit commutatif

$$
\CD
{\Cal C}_0(q_2) @> id     >> {\Cal C}_0(q_2)  \\
@V PE_{0,a,m'}  VV        @VV EP_{0,a,m}    V \\
{\Cal C}'_1(m') @> \psi^* >> {\Cal C}_1(m)
\endCD 
$$

\noindent Si tel est le cas, les invariants de tout couple 1-algébrisable induit $({\Cal X},F)_a$, sont majorés par ceux du couple 1'-algébrisable induit $(\chi',f')_a$ tel que 

$$
({\Cal X},F)_a=\psi^*(\chi',f')_a
$$

\noindent Or, ce dernier couple est l'éclaté d'un couple 0-algébrisable, dont les invariants peuvent être majorés, en utilisant la théorie de Khovanski [K], et les idées développées dans les travaux [IY], [NY].\par

\vskip 1cm

\noindent{\bf Remerciements.} Je tiens à remercier L. Gavrilov qui, par de brèves et sincères discussions, a stimulé mon intérêt pour le sujet des intégrales Abéliennes.\par

\vskip 1cm

\noindent{\bf R\'ef\'erences.}\par

\roster

\item"[AL]" Andronov A.A., Leontovitch E.A., Gordon I.I. and Maier A.G., {\sl Theory of bifurcations of dynamic systems on a plane}, (New York: Wiley) (1973), 1-482.\par

\item"[AV]" Arnold V.I., Varchenko A.N. and Gussein-Zade S.M., {\sl Singularities of differentiable mappings II, Monodromy and asymptotics of integrals}, Monographs in Mathematics vol. {\bf 82} (Boston: Birkhauser).\par

\item"[B]" Bonnet P., {\sl Description of the module of relatively exact 1-forms modulo a polynomial} $f$ {\sl on} $\CC^2$, Université de Bourgogne, IMB, Préprint 184 (1999).\par

\item"[G1]" Gavrilov L., {\sl Petrov modules and zeros of Abelian integrals}, Bull. Sci. Maths. {\bf 122} (1998), 571-584.\par

\item"[G2]" Gavrilov L., {\sl Abelian integrals related to Morse polynomials and perturbations of plane Hamiltonian vector fields}, Ann. Inst. Fourier {\bf 2} (1999).\par

\item"[G3]" Gavrilov L., {\sl The infinitesimal 16th Hilbert problem in the quadratic case}, Inv. Maths. {\bf 143} (2001), 449-497.\par

\item"[G4]" Gavrilov L., {\sl Higher order Poincaré-Pontryagin functions and iterated path integrals}, Ann. Fac. Sci. de Toulouse {\bf 14} (2005), 677-696.\par

\item"[I1]" Il'Yashenko J.S., {\sl The multiplicity of limit cycles arising from perturbations of the form $w'=P/Q$ of a Hamilton equation in the real and complex domain}, Trudy Sem. Petrovsk. {\bf 3} (1978), 49-60 (Engl. transl. Am. Math. Soc. Transl. {\bf 118} 191-202).\par

\item"[I2]" Il'Yashenko J.S., {\sl Appearance of limit cycles in perturbation of the equation $dw/dz=-R_z/R_w$ where $R(z,w)$ is a polynomial}, USSR Math. Sbornik {\bf 78} (1969), 360-373.\par

\item"[I3]" Il'Yashenko J.S., {\sl An exemple of equation $dw/dz=P_n(z,w)/Q_n(z,w)$ having a countable number of limit cycles and an arbitrary large genre after Petrovski-Landis}, USSR Math. Sbornik {\bf 80} (1969), 388-404.\par

\item"[IY]" Il'Yashenko J.S. and Yakovenko S., {\sl Counting real zeros of analytic functions satisfying linear ordinary differential equations}, J. Diff. Equations {\bf 126} (1995), no. 1, 87-105.\par

\item"[K]" Khovanski A.G., {\sl Fewnomials}, AMS Publi., Providence, RI, (1991).\par

\item"[Ma]" Mardesic P., {\sl An explicit bound for the multiplicity of zeros of generic Abelian integrals}, Nonlinearity {\bf 4} (1991), 845-852.\par 

\item"[Mo]" Mourtada A., {\sl Action de dérivations irréductibles sur les algèbres quasi-régulières d'Hilbert}, En cours de révision pour publication.\par

\item"[NY]" Novikov D. and Yakovenko S., {\sl Simple exponential estimate for the number of real zeros of complete Abelian integrals}, Ann. Inst. Fourier (Grenoble) {\bf 45} (1995), no. 4, 897-927.\par

\item"[R]" Rousseau C., {\sl Bifurcation methods in polynomial systems , Bifurcation and Periodic Orbits of Vector Fields}, ed. D. Schlomiuk (Dordrecht: Kluwer) (1993), 383-428.\par

\item"[Ro]" Roussarie R., {\sl Cyclicité finie des lacets et des points cuspidaux}, Nonlinearity {\bf 2} (1989), 73-117.\par

\item"[V]" Varchenko A.N., {\sl Estimate of the number of zeros of Abelian integrals depending on a parameter and limit cycles}, Funct. Anal. {\bf 18} (1984), 98-108.\par 

\item"[Y]" Yakovenko S., {\sl Complete Abelian integrals as rational envelopes}, Nonlinearity {\bf 7} (1994), 1237-1250.\par

\endroster

\enddocument